\documentclass[10pt, reqno]{amsart}

\usepackage{amssymb,latexsym}
\usepackage{eucal}
\usepackage{tikz}

\usepackage{amssymb,latexsym}
\usepackage{eucal}
\usepackage{tikz}
\usepackage{comment}

\newtheorem*{th1}{Theorem}
\newtheorem*{lm1}{Lemma 1}
\newtheorem*{lm2}{Lemma 2}

\newtheorem*{cr1}{Corollary 1}
\newtheorem*{cr2}{Corollary 2}
\newtheorem*{cr3}{Corollary 3}
\newtheorem*{cr4}{Corollary 4}
\newtheorem*{cr5}{Corollary 5}

\newcommand{\al}{\alpha}

\newcommand{\p}{\partial}

\newcommand{\ph}{\varphi}
\newcommand{\D}{\Delta}

\newcommand{\wtl}{\widetilde}

\DeclareMathOperator{\mg}{\textnormal{mg}}
\DeclareMathOperator{\eg}{\textnormal{eg}}
\DeclareMathOperator{\qp}{\textnormal{qp}}

\newcommand{\FF}{\mathcal{F}}

\newcommand{\bpp}{boundary path}
\newcommand{\cncm}{connected component}

\begin{document}

\title[Quasiperiodic and mixed commutator factorizations in free products]
{Quasiperiodic and mixed commutator factorizations in free products of groups}

 \author{Sergei V. Ivanov}
\address{Department of Mathematics\\
  University of Illinois \\
  Urbana\\   IL 61801\\ U.S.A.} \email{ivanov@illinois.edu}

  \author{Anton A. Klyachko}
\address{Department of Mechanics and Mathematics\\
  Moscow State University \\
 Moscow\\   119991 \\ Russia} \email{klyachko@mech.math.msu.su}
\thanks{The work of the second author was supported by the Russian Foundation for Basic Research, project  15-01-05823.}

\begin{abstract} It is well known that a nontrivial  commutator in a free group
is never a proper power. We prove a theorem that generalizes this fact and has
several worthwhile corollaries. For example,  an equation $[ x_1, y_1] \ldots [ x_k, y_k] = z^n$,  where $n \ge 2k$, in the free product  $\mathcal{F}$ of groups without nontrivial elements of order $\le n$ implies that $z$ is conjugate to an element of a free factor of $\mathcal{F}$. If a nontrivial commutator in a free group
factors into a product of elements which are conjugate to each other then all these
elements are distinct.
\end{abstract}

\maketitle

\section{Introduction}

It was observed by  Sch\"utzenberger \cite{Sch59} that  a  nontrivial commutator in a free
group is never a proper power. This result was generalized in different
directions: for values of other than commutator words  on free groups by Baumslag and Steinberg \cite{BS},
for values of commutators on free products of groups by Comerford, Edmunds, and Rosenberger \cite{CER94},
and for values of commutators on small cancellation groups by Frenkel and the second author \cite{FK12}.
Our Theorem could be considered as one more such a generalization.

\begin{th1} \label{th1} Let $G_\al$, $\alpha \in I$, be torsion free groups and  let
$\FF= \underset{\al \in I} * G_\alpha$ denote their free product. Suppose that
$$
c_1\dots c_k d_1\dots d_\ell =  h_1^{n_1}\dots h_m^{n_m}
$$
in $\FF$, where $c_1,\dots, c_k, d_1,\dots, d_\ell,  h_1,\dots, h_m$ are elements of $\FF$ such that
$c_i$ are commutators, $d_{i'}$ are conjugate to elements of free factors  of $\FF$,
$h_j$ are conjugate to each other and are not conjugate to elements of  free factors  of $\FF$, and
$n_j$ are positive integers. Then
$$
\sum_{j=1}^m{(n_j-1)} \le 2k+\ell-2 .
$$
In addition, the same statement holds for any free product $ \FF$ of groups
with torsion  whenever the order of
every letter of a cyclically reduced word conjugate to  $h_1$
is greater than $\sum_{j=1}^m{n_j}$.
\end{th1}

M. Culler \cite{Cull81} discovered that, in the free group $F(a,b)$ with free generators $a, b$, the cube $[a,b]^3$ of
the commutator $[a,b] := a^{-1}b^{-1}ab$ is a product of two commutators,
$$
[a,b]^3=[a^{-1}ba,a^{-2}bab^{-1}][bab^{-1},b^2] .
$$
Moreover, $[a,b]^n$ is a product of $k$ commutators whenever $n\le 2k-1$, see \cite{Cull81}.

Comerford, Comerford, and Edmunds  \cite{CCE91}  proved that a nontrivial product of two commutators
in a free group cannot be more than a third power, i.e., the
equality $[x_1,y_1] [x_2,y_2]=z^n$, where $n \ge 4$, in a free group implies that $z=1$.

The authors of \cite{CCE91} conjectured that, in a free group, the Culler's
bound $n\le 2k-1$ is sharp. In other words, the
Comerford-Comerford-Edmunds conjecture asserts that, in a free group, the
equality $[x_1,y_1]\ldots[x_k,y_k]=z^n$, where $n\ge 2k$, implies that
$z=1$.
This conjecture was proven by Duncan and Howie  \cite[Theorem 3.3]{DH91} by establishing that,
in the free product $A*B$ of two locally indicable groups $A, B$, the equality
$[x_1,y_1]\ldots[x_k,y_k]=z^n$, where $n\ge2k$, implies that $z$ is
conjugate to an element of $A$ or $B$.
Our Theorem implies, in particular, that a similar result holds true for the free
product of  torsion free groups and for the free
product  of groups with torsion bounded below. Note that, in the case of small torsion,
a similar result is no longer true: in the infinite dihedral group
$\langle c \rangle_2 * \langle d \rangle_2$, every power of the commutator $[c,d]$ is also a commutator.

\begin{cr1}\label{cr1}
Let $\FF = \underset{\al \in I} * G_\alpha$  be  the free product
of groups $G_\al$,  $\alpha \in I$, that contain no nontrivial elements
of order $\le n$. Then an equality
$[x_1, y_1] \ldots [x_k, y_k] = z^n$,  where $n \ge 2k$,
in $\FF$ implies that $z$ is conjugate to an element of a free factor of
$\FF$.
\end{cr1}

We remark that  Corollary~\ref{cr1} is similar to
a recent result of Chen \cite[Corollary~3.7]{C}, however,
neither Corollary~\ref{cr1} follows  from  \cite[Corollary~3.7]{C} nor
\cite[Corollary~3.7]{C} follows from Corollary~\ref{cr1}.

We now state more corollaries of our Theorem.

\begin{cr2}\label{cr2}
Let $\FF$ be the free product of torsion free groups $G_\al$,
$\alpha \in I$. Suppose that
$$
c_1\dots c_k d_1 \dots d_\ell= h_1\dots h_m,
$$
in $\FF$, where $c_1,\dots, c_k, d_1,\dots, d_\ell,  h_1,\dots, h_m$ are
elements of $\FF$ such that
$c_i$ are commutators, $d_i$ are conjugate to elements of free factors
of $\FF$,
and $h_i$ are conjugate to each other and are not conjugate to elements of
free factors  of $\FF$. Then
no element occurs  in the sequence $h_1,\dots,h_m$  more than
$2k+\ell-1$ times.
\end{cr2}

\begin{cr3}\label{cr3}
Suppose  $F$ is a free group and $c =  h_1 \dots h_m \ne 1$ in $F$, where
$c, h_1, \dots,h_m \in F$ are
such  that $c$ is a commutator and $h_1,\dots,h_m$ are
conjugate to each other.
Then  $h_1,\dots, h_m$ are all distinct.

More generally, suppose that $c_1\dots c_k  =  h_1 \dots h_m \ne 1$
in a free group $F$,
where $c_1, \dots,c_k,   h_1, \dots,h_m \in F$ are such that
$c_1, \dots,c_k$ are  commutators and $h_1,\dots,h_m$ are conjugate to
each other.
Then  no element  occurs  in the sequence $h_1,\dots,h_m$
more than $2k-1$ times.
\end{cr3}

\begin{cr4}\label{cr4}
Let $A*B$ be the free product of torsion free groups $A, B$. Then
no nontrivial element $a  \in A$ is a product of elements that are
conjugate to each other and are  not conjugate  to
an element of $A$.

More generally,  suppose that
$a_1 b_1\dots a_\ell b_\ell  =  h_1 \dots h_m \ne 1$ in  $A*B$,
where $a_1, \dots, a_\ell \in A \setminus \{ 1 \}$,
$b_1, \dots, b_\ell \in B\setminus \{ 1 \}$,
and $h_1,\dots,h_m \in A*B$ are conjugate to each other
and are not conjugate to an element of $A \cup B$. Then  no element
occurs  in the sequence $h_1,\dots,h_m$ more than $2\ell-1$ times.
\end{cr4}

Let $w$ be an element of the free product $\FF$ of groups $G_\al$,  $\alpha \in I$.
A {\em mixed commutator factorization} for $w$ is an equality in $\FF$ of the form
\begin{equation}\label{mcf}
  w = c_1 \dots c_k d_1  \dots d_\ell  ,
\end{equation}
where $c_i$ are commutators and $d_j$ are conjugate to elements of free factors  of $\FF$.
The {\em mixed genus}  $\mg(w)$ of  $w$ is defined to be a  minimal integer $s$ such that  $s= 2k+\ell$  over all mixed commutator factorizations \eqref{mcf} for $w$.

For example, if $\mg(w) \le 1$ then $w$ is conjugate to an element of a free factor of $\FF$ and if
$\mg(w) = 2$ then $w$ is a commutator or a product of two elements  conjugate to nontrivial elements of free factors  of $\FF$.
\smallskip

We remark that Culler \cite{Cull81} introduced the genus $g(w)$ for an element $w$ of the free product $A*B$ of two groups $A, B$
as a minimal number of commutators
needed to write $w$ as the product of these  commutators or $g(w) := \infty$ if $w$ is not a product of commutators.
Culler \cite{Cull81}  gave an algorithm that computes the genus $g(w)$
of $w$ whenever the genera of elements can be computed in free factors $A, B$. The genus $g(w)$ can be defined in the same fashion for an element $w$ of an arbitrary free product $\FF$ of groups.

Let a free group $F$ be  considered  as the free product of its cyclic subgroups.
Grigorchuk and Kurchanov \cite{GK} defined the width $h(w)$
of an element $w $  of $F$
as a minimal number of elements that are conjugate to elements of free factors of $F$ and that are needed to write $w$ as their product. Grigorchuk and Kurchanov \cite{GK}  gave an algorithm that computes the width $h(w)$ of $w \in F$, see also  \cite{Iv}, \cite{Ol}.  The width $h(w)$   can be defined in the same manner for an element $w$ of an arbitrary free product $\FF$ of groups.

It is worthwhile to note that our definition of the mixed genus $\mg(w)$ of an element $w$ of an arbitrary free product $\FF$  combines the foregoing two definitions and the number
$\mg(w)$ satisfies the inequalities $\mg(w) \le 2g(w)$ and $\mg(w) \le h(w)$.
However, it is not clear how to algorithmically compute the mixed genus   $\mg(w)$  even for elements of a free group.
\medskip

A {\em quasiperiodic   factorization} for an element $w$ of the free product  $\FF$ of groups $G_\al$,  $\alpha \in I$,  is an equality in $\FF$ of the form
\begin{equation}\label{psf}
    w = h_1^{n_1} \ldots h_m^{n_m} ,
\end{equation}
where
$h_1, \ldots, h_m$ are conjugate to each other and are not conjugate to an element of a free factor  $G_\al$,
$n_1, \ldots, n_m$ are positive  integers, 
and $m \ge 1$. 


The {\em quasiperiodicity}  $\qp(w)$ of $w$ is  defined to be a maximal integer $r$ such that $r = \sum_{j=1}^m(n_j-1)$ over all quasiperiodic   factorizations \eqref{psf} for $w$ if there are such factorizations and the set of such $r$ is bounded above. If the set of such $r$ is not bounded above, we set  $\qp(w) := +\infty$ and if there are no such factorizations for $w$, we set $\qp(w) := -\infty$.

It is clear that, for every $w \in \FF$ such that $w$  is not conjugate to an element of a free factor, we have
$\qp(w) \ge 0$. As another example, consider two elements $u, v \in \FF$ that are conjugate and are not conjugate to  an element of a free factor of $\FF$. Then $\qp( u^4 v^2) \ge 4$ and  $\qp( u^3  v u v) \ge 3$ as $u^3  v u v = u^4 v^u v $, where $v^u := u^{-1} v u$.

Note that if $w = h_1 \ldots h_m$ in $\FF$, where $h_1, \ldots, h_m$ are conjugate to each other and are not conjugate to an element of a free factor  $G_\al$, and $s$ elements among   $h_1, \ldots, h_m$ are equal each other, then $\qp( w) \ge s-1$. Indeed, we can apply the identity $uv=vu^v$ and rearrange the
factors  $h_1, \ldots, h_m$  in such a way that the equal $s$  factors
would form an $s$th power. This observation,  in particular, implies that,
if the free product $A*B$ has torsion, then   $\qp(1) = +\infty$. Indeed, if
an element $a\in A$ has order $m>1$ and $b \in B$ is nontrivial then
$$
1=[a,b][a,b]^{a^{-1}} \ldots [a,b]^{a^{1-m}} =
([a,b][a,b]^{a^{-1}}  \ldots [a,b]^{a^{1-m}})^{2016} .
$$

These equalities mean that  $\qp(1) = \qp([a,b]) =+\infty$.
(It is not clear what could be $\qp(a), \qp(ab)$ in this situation.)
On the other hand, for free products of groups without torsion we have a
nicer situation.

\begin{cr5}\label{cr5} Let $\FF$ be the free product of torsion free groups  $G_\al$,  $\alpha \in I$. Then, for every
$w \in \FF$,  the quasiperiodicity  $\qp(w)$ of $w$ satisfies $\qp(w) \le \mg(w)-2 <  +\infty$.
Furthermore,  
$\qp(w) = -\infty$ if and only if $w$ is  conjugate to an element of a free factor of $\FF$, otherwise, $\qp(w) \ge 0$  is finite.  
\end{cr5}

We remark that the bound $\qp(w) \le \mg(w)-2$ of Corollary~5 is sharp
as follows from the equality
$$
(ab)^n=a^nb^{a^{n-1}}b^{a^{n-2}}\dots b^ab
$$
that proves that if $a \in A$, $b \in B$  are nontrivial then
$\qp( (ab)^n ) \ge n-1$ and $\mg((ab)^n ) \le n+1$. The sharpness of the bound $\qp(w) \le \mg(w)-2$  also follows from
the Culler's observation \cite{Cull81}
that $[a,b]^n$ is a product of $k$ commutators whenever $n \le 2k-1$.

Our arguments utilize  diagrams over free products of groups and are based on a car-crash lemma of \cite{Kl93},  \cite{Kl97}, \cite{Kl05},  see also \cite{FeR96}, that
has had quite a few applications in group theory, see   \cite{CG95},  \cite{CG00}, \cite{CR01},  \cite{FeR98},  \cite{FK12}, \cite{FoR05},  \cite{IK},  \cite{Kl06a}, \cite{Kl06b},
\cite{Kl07}, \cite{Kl09}, \cite{KlL12}, \cite{Le09}.

In Sect. 2, we define  diagrams over free products of groups and prove
a lemma on geometric meaning of the mixed genus.
In Sect. 3, we state a car-crash lemma.
Sect. 4 contains the  proof of our Theorem.

\section{Preliminaries}

Suppose that $S$ is an oriented compact  closed
surface. Note that $S$ need not be connected.

A {\em map} on $S$ is a finite 2-complex $\D$ embedded into $S$.
We call $S$ the {\em underlying surface} for  $\D$,  denoted $S = S(\D)$.
If the embedding of $\D$ into $S$ is surjective, i.e., $\D$ has no boundary,  we say that the map $\D$ is {\em closed}.

The set of $i$-cells of a finite 2-complex $\D$ is denoted $\D(i)$, $i=0,1,2$. The closures of
$i$-cells  of $\D$  for $i=0,1,2$ are called {\em vertices, edges, faces}, resp. The 1-skeleton of $\D$, consisting of vertices and edges, is a graph denoted $\D[1]$.

If $F$ is a face of a map $\D$ then a boundary path $\p F$ of $F$
is oriented in positive, i.e., in counterclockwise,
direction. Recall that $S(\D)$ is oriented.
If $\p F = e_1e_2 \ldots e_k$, where $e_1,e_2, \ldots, e_k$
are oriented edges, then the subpaths
$e_1e_2, e_2e_3, \ldots, e_k e_1$
of  $\p F$  are called  {\em corners}
of $F$.
If $e_i e_{i+1}$ is a corner of $F$ then
the terminal vertex of $e_i$ is called the {\em vertex} of
$e_i e_{i+1}$ and is denoted $\nu(e_ie_{i+1})$.

If $e$ is an oriented  edge of a 2-complex $\D$ then $e_-, e_+$
denote the initial, terminal, resp., vertices of $e$.
By $e^{-1}$ we mean the edge with opposite to $e$ orientation.
If $p = e_1 \dots e_k$ is a path in $\D$, where
$e_1, \ldots, e_k$ are oriented edges, then the initial and
terminal vertices of $p$ are defined by  $p_- := (e_1)_-$ and
$p_+ := (e_k)_+$, resp., and  $p^{-1} :=  e_k^{-1} \ldots e_1^{-1}$.

Let $C(\D)$ denote the set of all corners of faces of a map $\D$ and let $A *B$ be the free product of two nontrivial groups $A, B$, where $A \cap B = \{ 1\}$.
A map $\D$ is called a {\em  diagram} over  $A*B$ if   $\D$ is equipped with two labeling functions
$$
\ph : C(\D) \to A \cup B, \qquad \theta : \D(0) \to \{ A, B\}
$$
and the following conditions are satisfied.

\begin{enumerate}
\item[(D1)] If $u, v$ are two vertices of $\D$ connected by an edge then
$\theta(u) \ne \theta(v)$. In particular,  $\D[1]$ is a  biparite graph.

\item[(D2)] For every corner $ee' \in C(\D)$ such that  $\theta( \nu(ee')) =A$, we have $\ph(ee') \in A$ and,
for every corner $ee' \in C(\D)$  such that  $\theta( \nu(ee')) =B$, we have $\ph(ee') \in B$.
\end{enumerate}

We remark that our definition of a diagram over  $A*B$ is different from the definitions of diagrams
over free products of groups used in books \cite{LS}, \cite{Ol89} and is similar to the definition introduced in
Howie's articles \cite{How83}, \cite{How90}.
\smallskip

Let $F$ be a face of a diagram $\D$ over  $A*B$  and let
$\p F = e_1e_2 \ldots e_k$, where $e_1,e_2, \ldots, e_k$
are oriented edges, be a boundary path of $F$.
A {\em label}  $\ph(\p F)$  of  $F$ is defined by setting
$$
\ph(\p F) :=  \ph( e_1 e_2) \ph( e_2 e_3) \ldots \ph(e_k e_1) ,
$$
i.e.,  $\ph(\p F)$ is the  product of consecutive, in positive direction, $\ph$-labels of corners of $F$.
It is clear that $\ph(\p F)$ is a word over  the alphabet $A \cup B$ and
$\ph(\p F)$ is defined up to a cyclic permutation.

Let $p = e_i e_{i+1} \ldots e_{i+\ell}$ be a subpath of a boundary path $\p F$ of a face $F$, where indices are modulo $k = |\p F|$. We define the  {\em label} $\ph(p)$ of $p$ to be the word
$$
\ph(p) :=  \ph(e_i e_{i+1} ) \ph(e_{i+1} e_{i+2}) \ldots \ph( e_{i+\ell-1} e_{i+\ell}) .
$$

If $v \in \D(0)$ is a vertex in the interior of a diagram $\D$ over $A*B$,
i.e., $v \not\in \p \D$, then a {\em label} $\ph(v)$ of $v$ is
the product  of $\ph$-labels of consecutive, in negative direction,  corners whose vertex is $v$.
We say that $v$ is an {\em $A$-vertex} if $\theta(v) =A$ and $v$ is a {\em $B$-vertex} if $\theta(v) =B$.
It is clear from the definitions that  $\ph(v) \in A$ if $v$ is an $A$-vertex and  $\ph(v) \in B$
if $v$ is a $B$-vertex. It is also clear that $\ph(v)$ is defined up to conjugation in $A$ or $B$.
If $\ph(v)$ is defined and $\ph(v) =1$ in $A$ or in $B$, depending on type of $v$, then we say that $v$ is a {\em regular} vertex. If $\ph(v)$ is defined and $\ph(v) \ne 1$ in $A$ or $B$ then
   $v$ is called an {\em irregular} vertex.  Note that a label $\ph(v)$ is not defined
for a vertex $v$ on the boundary $\p \D$ of $\D$.

We remark that similar  diagrams were considered in
\cite{How83}, \cite{How90}, \cite{Kl93}, \cite{Le09} and some other papers but our
definitions are slightly different.

For example, the diagram  depicted in Fig.~1 has a torus as the  underlying surface and it
is drawn as a rectangle with  opposite sides to be identified. This diagram contains two vertices, three
edges, one face, and three corners with $\ph$-label $a \in A$ and three corners with $\ph$-label $b\in B$.
If the vertices are regular, then $a^3 =1 $ in $A$
and $b^3=1$ in $B$. The label of the face is
$(ab)^3$. This diagram demonstrates that if $a \in A$ and $b \in B$ have order 3 then
$(ab)^3$ is a commutator. A complete description of commutators in a free product
of groups that are not conjugate to elements of free factors and are proper
powers is given in \cite{CER94}.
\medskip

\begin{center}
\begin{tikzpicture}[scale=.8]
\draw  (-5,3)[dashed] rectangle (-.5,1);
\draw (-5,2) -- (-0.5,2);
\draw  (-3.5,2) [fill] circle (0.05);
\draw  (-1.5,2)[fill] circle (0.05);
\draw (-3.5,2) -- (-2.5,3);
\draw (-2.5,1) -- (-1.5,2);
\node at (-2.5,0.2) {Fig. 1};
\node at (-3.7,2.3) {$a$};
\node at (-2.9,2.3) {$a$};
\node at (-3.5,1.7) {$a$};
\node at (-1.2,1.7) {$b$};
\node at (-1.5,2.3) {$b$};
\node at (-2.2,1.7) {$b$};
\end{tikzpicture}
\end{center}

We call a  diagram $\D$  over $A*B$ \emph{reduced} if  $\D$ has no face with a corner whose
$\ph$-label is 1.

The \emph{extended genus} $\eg(\D)$    of a  diagram $\D$ over  $A*B$   is defined by
\begin{equation}\label{eg}
 \eg(\D) :=  2-\chi(\D) + r_0(\D) ,
\end{equation}
where  $\chi(\D) = |\D(0)| - |\D(1)| + |\D(2)| $ is  the Euler
characteristic of $\D$ and $r_0(\D)$ is the number of irregular vertices in $\D$.
\medskip

We consider  elements of the free product $A *B$  as words over the alphabet $A \cup B$,
where $A \cap B = \{ 1\}$,  whose elements are called letters.
A word $w = a_1 \ldots a_\ell$, where $a_1, \ldots, a_\ell \in A\cup B$ are letters, is called \emph{reduced} if
$w$ is nonempty, none of the letters $a_1, \ldots, a_\ell$ is 1 and, for every $i$, the letters $a_i, a_{i+1}$ do not belong to the same free factor of $A*B$. The \emph{length} of a word $w$ is the number of letters in $w$, denoted  $|w|$.
A word $w$ is  \emph{cyclically reduced} if $w$ is nonempty and $w^2$ is reduced.

The definitions of  reduced  and cyclically reduced words
that represent elements of an arbitrary  free product   $\FF = \underset {\al \in I} \ast G_\al$
of nontrivial groups $G_\al$  are analogous.

If $u, w$ are two words over $A \cup B$, then $u \equiv w$ means the literal or letter-by-letter equality of words.
If words $u, w$  are equal as elements of $A *B$, we write $u \overset {*} = w$.

If $w$ is a word over  $A \cup B$, we let $\delta_1(w)$ denote the word obtained from $w$ by deletion of all occurrences of the  letter $1 \in A \cup B$. By writing $w \equiv_1  u$ we mean that $\delta_1(w) \equiv \delta_1(u)$. If $a, b \in A$ then by writing
$a \overset A = b$ we mean that $a = b$ in $A$. Similarly, the notation $a \overset B = b$ means that
$a, b \in B$  and $a = b$ in $B$.
\medskip

We mention without proof that the mixed genus  $\mg(w)$ of a cyclically
reduced word  $w$ over $A \cup B$  is equal to the
minimal  extended genus $\eg(\D)$  of a reduced closed diagram $\D$ over $A*B$
that contains a single face whose label is the word  $w$.
In  this paper we need only the inequality $\mg(w) \ge \eg(\D)$ that
follows from Lemma~1 below. We also remark that the inequality $\mg(w) \ge \eg(\D)$ of Lemma~1 is actually an equality but we will not need this equality either.

\begin{lm1}\label{lem1}
Suppose that $u_1, \ldots, u_m$ are nonempty cyclically reduced
words over the alphabet $A \cup B$. Then the minimal
mixed genus $\mg(w)$ of an element  $w \in A*B$ of the form
$w \overset {*} = s_1 u_1 s_1^{-1} \ldots s_m  u_m s_m^{-1}$
satisfies  $\mg(w) \ge \eg(\D)$, where $\eg(\D)$ is  the minimal
extended genus $\eg(\D)$ of a  reduced closed diagram $\D$ over $A*B$
that contains precisely $m$ faces whose labels are the
words $u_1, \ldots, u_m$.
\end{lm1}

\begin{proof}  Suppose that $w_0$ is a word of a minimal mixed genus
$n = \mg(w_0)$ among all words $w$ of the form
$$
w \equiv s_1 u_1 s_1^{-1} \dots s_m u_m s_m^{-1} ,
$$
where each $s_i \in A*B$ is a reduced word or $s_i \equiv 1$. Since $n = \mg(w_0)$, there is a factorization for $w_0$ of the form
$$
w_0 \overset {*} =  [v_1, t_1]  \dots [v_k, t_k]   d_1  \dots d_\ell  ,
$$
where $n = 2k +\ell$, for every $i$, $v_i, t_i$ are  reduced words, and, for each $j$,
$d_j \equiv d_{j,1} d_{j,0}d_{j,1}^{-1}$, $d_{j,0} \in A \cup B$, $d_{j,0} \ne 1$,
$d_{j,1}$ is a reduced word or $d_{j,1}  \equiv  1$.

Consider the word
\begin{equation}\label{eq0}
s_m u_m^{-1} s_m^{-1} \dots   s_1 u_1^{-1} s_1^{-1} v_1^{-1} t_1^{-1} v_1 t_1 \dots  v_k^{-1} t_k^{-1} v_k t_k d_1 \dots d_\ell .
\end{equation}

Let $\D_0$ be a diagram  over $A*B$ that consists of a single face $H$ whose boundary path $\p H$ has the following factorization
$$
\p H = p_1 q_1 p_2 q_2 \ldots p_L q_L ,
$$
where $L = 3m +4k +\ell$, $p_1, q_1, \ldots p_L q_L$ are subpaths of  $\p H$,
and the sequence of words $\ph(p_1), \ldots, \ph(p_L)$ is identical to the sequence
of subwords  $s_1, u_1^{-1}, s_1^{-1}, \dots,  d_\ell$  distinguished in the word \eqref{eq0}.
The paths $q_1, \ldots, q_L$ have labels equal to powers of the letter $1 \in A \cup B$ and $|q_i| = 2$ or $|q_i| = 3$,  hence,
$\ph(q_i) \equiv 1^{|q_i|-1}$ with $|q_i|-1 \ge 1$. The $\ph$-label of a corner of $H$ whose vertex is $(p_i)_-$ or
$(p_i)_+$ is also $1$. It is easy to see that we can assign $\theta$-labels to corners of $H$ so that both properties
(D1)--(D2) hold true and $\D_0$ is indeed a  diagram  over $A*B$. Note that the choice between  $|q_i| = 2$ or $|q_i| = 3$ depends on  $\theta((p_i)_+)$ and  $\theta((p_{i+1})_-)$, here indices are modulo $L$.

Let $x$ be one of the subwords $s_m, u_m^{-1}, s_m^{-1},  \dots,  d_\ell$  distinguished in the word \eqref{eq0} and let
$p(x)$ denote the corresponding path among  $p_1 , p_2, \ldots, p_L$ such that  $\ph(p(x)) \equiv x$, i.e.,  we assume that
$p(s_m) = p_1$, $p(u_m^{-1}) = p_2, \dots$,  $p(d_\ell) = p_L$.
\medskip

We now make some surgeries over $\D_0$. We remark that $\theta$-labels of vertices never change under these surgeries.

Observe that the subpath $p(d_j)$ of $\p H$ has even length because $d_j \equiv d_{j,1} d_{j,0}d_{j,1}^{-1}$,
where $d_{j,0} \in A \cup B$ and $d_{j,0} \ne 1$, $d_{j,1}$ is a reduced word or  $d_{j,1} \equiv  1$. Hence, there is a factorization $p(d_j) = p_1(d_j) p_2(d_j)$, where $| p_1(d_j)|  = | p_2(d_j)|$ and $\ph(p_1(d_j)) \equiv \ph(p_2(d_j))^{-1}$. Therefore, by identifying $p_1(d_j)$ and $p_2(d_j)^{-1}$ for each $j=1, \dots, \ell$,
within $H$, so that  the subpath $p(d_j)$  of $\p H$ turns into
 $\bar p_1(d_j) \bar p_1(d_j)^{-1} $, see Fig.~2,    we obtain a diagram $\D_1$  over $A*B$ with a single face, still denoted
$H$, and $\ell$ irregular vertices $(\bar p_1(d_j) )_+$, $j=1, \dots, \ell$. Note that all vertices of the paths
$\bar p_1(d_j)$, except for their end vertices $(\bar p_1(d_j))_- \in \p \D_1$ and  $(\bar p_1(d_j))_+$ are regular.

\begin{center}
\usetikzlibrary{arrows}
\begin{tikzpicture}[scale=.95]
\draw (-3.5,2) -- (-2.5,2.);
\draw[-latex] (-3.1,2) -- (-2.9,2);
\draw[-latex] (-1.45,0) -- (-1.4,0);
\draw[-latex] (-1.71,1.33) -- (-1.54,1.33);
\draw  (-3.5,2) [fill] circle (0.06);
\draw  (-2.5,2)[fill] circle (0.06);
\draw  (-1.5,4) [fill] circle (0.06);
\draw  (-1.5,3)[fill] circle (0.06);
\draw (-1.5,4) -- (-1.5,3);
\node at (-4.4,.4) {Fig. 2};
\node at (-1.7,1.65) {$\partial G_i$};
\node at (-2.9,2.4) {$\bar p(s_i)$};
\node at (-1.5,.4) {$\partial H$};
\node at (-1.7,2.2) {$G_i$};
\node at (-.2,2.) {$H$};
\node at (-2.24,2.) {$1$};
\node at (-1.97,3.) {$d_{j,0}$};
\node at (-.81,3.4) {$\bar p_1(d_j) $};
\draw[-latex] (-1.5,3.4) -- (-1.5,3.37);
\draw  (-1.5,2) ellipse (2 and 2);
\draw  (-1.7,2) ellipse (.8 and .67);
\end{tikzpicture}
\end{center}

Our next step is to identify, for every $i=1, \dots, m$, the path $p(s_i)$ with $p(s_i^{-1})^{-1}$ within $H$, getting thereby
a new path $\bar p(s_i)$ and a new map $\D_2$ with $m+1$ faces  $H, G_1, \dots, G_m$ such that the \bpp\
$\p |_{(\bar p(s_i))_+ } G_i$ of $G_i$ starting at the vertex $(\bar p(s_i))_+$ is a subpath of $(\p H)^{-1}$, see Fig.~2.
We assign 1 as the $\ph$-label to the corner of $G_i$ whose  vertex is $(\bar p(s_i))_+$ and, to every other corner of
$G_i$, we assign  $\ph$-label equal to the inverse of the $\ph$-label of the corner of $H$ with the same vertex.
Recall that $\theta$-labels of vertices do not change.
Such assignments produce a  diagram $\D_2$  over $A*B$
without additional irregular vertices because all vertices of
$\p G_i$ and $\bar p(s_i)$ are regular.

We now identify the path $p(u_i^{-1})$ with $p(u_i)^{-1}$ and the path $p(t_i^{-1})$ with $p(t_i)^{-1}$ for every
$i =1, \dots, k$.
Doing these identifications, results in a diagram $\D_3$  over $A*B$ on an oriented surface of genus $k$ such that $\D_3$ consists of  $m+1$ faces
$H, G_1, \dots, G_m$, $\D_3$ has a single boundary component, denoted  $\p \D_3$, and $\D_3$
contains $\ell$ irregular vertices.  The images of paths
$p(u_i)$,  $p(t_i)$ in $\D_3$ are denoted  $\bar p(u_i)$,  $\bar p(t_i)$, resp., $i   =1, \dots, k$.

Note that the vertices of paths $\bar p(u_i)$,  $\bar p(t_i)$, different from their end vertices, are all regular and
the end vertices of $\bar p(u_i)$,  $\bar p(t_i)$ belong to the boundary path $\p \D_3$.
We also observe that if $ee' \in C(\D_3)$ is a corner whose vertex belongs to  $\p \D_3$ then
 $\ph(e e') = 1$. Therefore, we may attach a new face $G_0$ such that
$|\p G_0| = |\p \D_3 |$ and  $  \ph(\p G_0) \equiv 1^{|\p \D_3 |}$ to $\p \D_3$ by identifying the paths $\p G$ and $\p \D_3$.

This attachment of $G_0$ to $\D_3$ produces a new diagram $\D_4$  over $A*B$ such that $\D_4$ is  closed,  $\chi(\D_4) = 2-2k$, $\D_4$  has $\ell$ irregular vertices,
$$
\eg(\D_4) = 2k+\ell = \mg(w_0) ,
$$
$\D_4$ contains $m$ faces $G_1,  \dots, G_m$ such that
$\ph( \p  G_i ) \equiv_1 u_i$, where $i=1, \dots, m$, and $\D_4$ contains two more faces $H, G_0$ such that
$\ph( \p  H ) \overset{*} =\ph( \p G_0 ) \overset{*} = 1$.

Thus we have constructed a closed diagram $\D_4$ over $A*B$ with some desired properties except for the properties of being reduced and having precisely $m$ faces whose labels are the words $u_1, \dots, u_m$.
\smallskip

For a closed  diagram $\D$ over $A*B$  consider the parameter
$$
\tau(\D) := (-\chi(\D), r_0(\D),  |\D(1)| ) ,
$$
where as above $\chi(\D)$ is the Euler characteristic of $\D$,  $r_0(\D)$ is the number of irregular vertices of $\D$,
and $|\D(1)|$ is the number of nonoriented edges of $\D$.
We partially order diagrams $\D$  over $A*B$ according to their parameters
$\tau(\D)$ which are ordered lexicographically, i.e., $\tau(\D) <  \tau(\D')$ if and only if
$-\chi(\D) < -\chi(\D')$ or   $-\chi(\D) = -\chi(\D')$ and   $r_0(\D) < r_0(\D')$ or $-\chi(\D) = -\chi(\D')$ and   $r_0(\D) = r_0(\D')$  and $ |\D(1)| < |\D'(1)|$.

Initializing, we set $\wtl \D := \D_4$ and note that
$$
\tau(\wtl \D) = ( 2-2k, \ell, \wtl \D(1)) , \quad  \eg(\wtl \D) \le \mg(w_0) .
$$

In our inductive arguments below we  do not assume that $\wtl \D$ is necessarily connected but we do
assume that $\wtl \D$ has the following property.

\begin{enumerate}
  \item[(P)] Every \cncm\ of a diagram $\D$  over $A*B$ contains a face $F$  such that $\ph(\p F) \equiv_1 u_i$
  for some $i=1,\dots, m$.
\end{enumerate}

Note that the number of connected components of a diagram $\D$ over $A*B$ with property (P) is at most $m$. Hence,
$-\chi(\D) \ge -2m$ because $-\chi(\D) \ge -2$ whenever $\D$ is connected. Since the second and the third components
of $\tau(\D)$  are nonnegative integers, it follows that there is no strictly decreasing infinite chain
$$
\tau(\D_1) > \tau(\D_2) > \ldots
$$
in which diagrams $\D_1, \D_2, \ldots$ have property (P).  This means that we may use induction on
parameter $\tau(\D)$ in our arguments below if intermediate diagrams, similarly to  $\wtl \D$,  also have  property (P).
\smallskip

Now we will make more surgeries over  $\wtl \D$ aimed to get a reduced diagram.

If $\wtl \D$ is reduced and has property (P)  then $\wtl \D$ is a required diagram and our proof is complete.
\medskip

Suppose that there is a corner $ef$ of a face $F$ of $\wtl \D$ such that $\ph(ef)=1$.
Consider three possible cases.

\smallskip
Case 1: Assume that $e = f^{-1}$, i.e.,  the vertex $ e_+ =f_-$ has degree 1 and the corner $ef = ee^{-1}$
is the only corner in $\wtl \D$  whose vertex is $e_+$.

If the degree of the vertex $e_-$ is also 1 then  the \cncm\ of
$\wtl \D$  that contains $e, f$ is a sphere that contains  the single face $F$ such that $| \p F| = 2$ and $\ph( \p F) \equiv 1 c$, where $c $ is the $\ph$-label of the second corner of $F$. Since $u_1, \dots, u_m$ are cyclically reduced words, it follows that the label of $F$ may not be one of $u_1, \dots, u_m$. This contradiction to property (P) of $\wtl \D$  proves that the degree of $e_-$ is greater than 1.
Hence, we may take  the edges $e, f$ out of $\wtl \D$ creating thereby a diagram $\wtl \D_1$ with property (P) and
$\eg(\wtl \D_1) =  \eg(\wtl \D)$.
The two consecutive corners $e'e$, $e^{-1} f'$ of $F$ will disappear  and, in their place, we obtain a single corner $e'f'$ whose $\ph$-label is defined by
$\ph(e'f') := \ph(e'e)\ph(ff')$, see Fig.~3, where  $\ph(e'e) = a_1$, $\ph(ff') = a_2$  and $a_1, \dots, a_4 \in A \cup B$.
In view of inequality $\tau(\wtl \D_1) < \tau(\wtl \D)$ and  $\eg(\wtl \D_1) =  \eg(\wtl \D)$, we can use the induction hypothesis and Case 1 is complete.

\begin{center}
\usetikzlibrary{arrows}
\begin{tikzpicture}[scale=1]
\draw (-4,2) -- (-0.5,2);
\draw[-latex] (-3.3,2) -- (-3.4,2);
\draw  (-4,2) [fill] circle (0.06);
\draw  (-2,2)[fill] circle (0.06);
\draw (-2,2) -- (-1.2,3);
\draw  (-2,2)--(-1.2,1) ;
\node at (.8,.4) {Fig. 3};
\node at (-4.49,2) {$1$};
\node at (-2.2,1.6) {$a_1$};
\node at (-3.2,2.4) {$e=f^{-1}$};
\node at (-1.2,1.6) {$a_4$};
\node at (-2.2,2.4) {$a_2$};
\node at (-1.2,2.4) {$a_3$};
\draw [thick][-latex](.5,2) -- (1.,2);
\draw  (3,2)[fill] circle (0.06);
\draw (3,2) -- (3.8,3);
\draw  (3,2)--(3.8,1) ;
\draw  (3,2)--(4.5,2) ;
\node at (2.38,2) {$a_1a_2$};
\node at (3.8,1.6) {$a_4$};
\node at (3.8,2.4) {$a_3$};
\end{tikzpicture}
\end{center}

\smallskip
Case 2: Suppose  $e_-\ne f_+$.

In this case we fold the edges $e$ and $f^{-1}$ within $F$, i.e., we identify $e$ and $f^{-1}$
through  the ``corner" of $F$ between them. The vertices $e_-, f_+$ become identical  and the two corners $e'e$, $ff'$ of $F$, whose vertices were $e_-, f_+$ before the fold, turn into a single corner $e'f'$ whose $\ph$-label is defined by
$\ph(e'f') := \ph(e'e)\ph(ff')$, see Fig.~4, where  $\ph(e'e) = a_1$, $\ph(ff') = a_5$  and $a_1, \dots, a_6 \in A \cup B$.
As a result, we obtain a  diagram $\wtl \D_1$  over $A*B$ such that
$\tau(\wtl \D_1) < \tau(\wtl \D)$ and $\eg(\wtl \D_1) =  \eg(\wtl \D)$.
By the induction hypothesis, Case 2 is  complete.

\begin{center}
\usetikzlibrary{arrows}
\begin{tikzpicture}[scale=.82]
\draw (-0.7,-.7) -- (1.,1.);
\draw (-.7,.7) -- (1.,-1.);
\draw (2.,-1.) -- (1.,-1.);
\draw (1,-2) -- (1.,-1.);
\draw (1,1) -- (2.,1.);
\draw (1,1) -- (1.,2);
\draw[-latex] (.5,.5) -- (.4,.4);
\draw[-latex] (.5,-.5) -- (.6,-.6);
\draw  (0,0) [fill] circle (0.06);
\draw  (5.5,0)[fill] circle (0.06);
\draw  (1,1) [fill] circle (0.06);
\draw  (1,-1)[fill] circle (0.06);

\node at (3.6,-2.4) {Fig. 4};
\node at (.6,0) {$1$};
\node at (1.3,.6) {$a_1$};
\node at (.2,.75) {$e$};
\node at (.2,-.75) {$f$};
\node at (.6,1.3) {$a_2$};
\node at (1.5,1.5) {$a_3$};
\node at (.6,-1.2) {$a_4$};
\node at (1.3,-.6) {$a_5$};
\node at (1.5,-1.5) {$a_6$};

\node at (6.7,.4) {$a_2$};
\node at (7.45,.7) {$a_3$};
\node at (8,0) {$a_1a_5$};
\node at (6.7,-.4) {$a_4$};
\node at (7.45,-.6) {$a_6$};

\draw [thick][-latex](3,0) -- (4,0);
\draw  (7.2,0)[fill] circle (0.06);
\draw (7.2,0) -- (7.,1);
\draw  (7.2,0)--(7.,-1) ;
\draw (7.2,0) -- (8.,.7);
\draw  (7.2,0)--(8.,-.7) ;
\draw  (5.5,0)--(7.2,0) ;
\draw (5.5,0) -- (4.8,.7);
\draw  (5.5,0)--(4.8,-.7) ;

\node at (1.7,0) {$H$};
\end{tikzpicture}
\end{center}

Case 3: Suppose   $e \ne  f^{-1}$ and $e_- = f_+$.

In this case, the path $ef$ is closed and defines a simple closed curve on the underlying surface $S(\wtl \D )$.  We cut $\wtl \D$ along this curve and obtain a new diagram $\wtl \D_0$ with two boundary components, oriented clockwise, which we denote by $e' f'$ and $(e'' f'')^{-1}$, where $e', e'' $ are the images  of $e$ in $\wtl \D_0$, $f', f'' $ are the images of $f$ in $\wtl \D_0$, and $e' f'$ is the image of the corner $ef$ of $H$ in  $\wtl \D_0$, see Fig.~5.

\begin{center}
\begin{tikzpicture}[scale=.67]
\draw (0,1) -- (.5,1.5);
\draw (0,1) -- (-.5,1.5);
\draw (0,-1) -- (.5,-1.5);
\draw (0,-1) -- (0,-1.7);
\draw (0,-1) -- (-.5,-1.5);
\draw (0,-1) -- (-.3,-.5);
\draw (0,-1) -- (.3,-.5);

\draw[-latex] (-1,-0.1) -- (-1,.1);
\draw[-latex] (1,0.1) -- (1,-.1);
\draw  (0,-1) [fill] circle (0.075);
\draw  (0,1) [fill] circle (0.075);
\node at (7,-2.7) {Fig. 5};
\node at (0,.60) {$1$};
\node at (-1.4,0) {$e$};
\node at (1.4,0) {$f$};

\draw [thick][-latex](2,0) -- (2.7,0);
\draw [thick][-latex](9.7,0) -- (10.3,0);

\node at (0,0) {$H$};
\draw (0,0)  circle (1);
\begin{scope}[xshift=4.5cm,yshift=0.4cm]
\draw (0,-1) -- (-.3,-.5);
\draw (0,-1) -- (.3,-.5);
\draw[-latex] (-1,-0.1) -- (-1,.1);
\draw[-latex] (1,0.1) -- (1,-.1);
\draw  (0,-1) [fill] circle (0.075);
\draw  (0,1) [fill] circle (0.075);
\node at (0,.60) {$1$};
\node at (-1.4,0) {$e'$};
\node at (1.4,0) {$f'$};
\node at (0,0) {$H$};
\draw (0,0)  circle (1);
\end{scope}

\begin{scope}[xshift=7.7cm,yshift=-0.4cm]
\draw (0,1) -- (.5,1.5);
\draw (0,1) -- (-.5,1.5);
\draw (0,-1) -- (.5,-1.5);
\draw (0,-1) -- (0,-1.7);
\draw (0,-1) -- (-.5,-1.5);
\draw[-latex] (-1,-0.1) -- (-1,.1);
\draw[-latex] (1,0.1) -- (1,-.1);
\draw  (0,-1) [fill] circle (0.075);
\draw  (0,1) [fill] circle (0.075);
\node at (-1.4,0) {$e''$};
\node at (1.4,0) {$f''$};
\draw (0,0)  circle (1);
\end{scope}

\begin{scope}[xshift=11.cm,yshift=0.4cm]
\draw (0,-1) -- (-.3,-1.5);
\draw (0,-1) -- (.3,-1.5);
\draw (0,-1) -- (0,1);
\draw[-latex] (0,0) -- (0,.1);
\draw  (0,-1) [fill] circle (0.075);
\draw  (0,1) [fill] circle (0.075);
\node at (0,1.5) {$1$};
\node at (1.46,0) {$e' = (f')^{-1}$};
\node at (-0.5,0.7) {$H$};
\end{scope}

\begin{scope}[xshift=14cm,yshift=-0.4cm]
\draw (0,-1) -- (0,1);
\draw[-latex] (0,0) -- (0,.1);
\draw  (0,-1) [fill] circle (0.075);
\draw  (0,1) [fill] circle (0.075);
\node at (1.57,0) {$e'' = (f'')^{-1}$};
\draw (0,1) -- (.5,1.5);
\draw (0,1) -- (-.5,1.5);
\draw (0,-1) -- (.5,-1.5);
\draw (0,-1) -- (0,-1.7);
\draw (0,-1) -- (-.5,-1.5);
\draw  (0,-1) [fill] circle (0.075);
\draw  (0,1) [fill] circle (0.075);
\end{scope}
\end{tikzpicture}
\end{center}

Note that
$\chi(\wtl \D_0) = \chi(\wtl \D)$ and the closed paths $e'f'$,  $e''f''$ might belong to different \cncm s of $\wtl \D_0$
which happens when $ef$ defines a separating curve on $S(\wtl \D )$.

We identify the edges $e'$ and $(f')^{-1}$ and  the edges  $e''$ and $(f'')^{-1}$ thus eliminating the boundary of
$\wtl \D_0$. The result is a closed diagram $\wtl \D_1$  over $A*B$ such that
$$
\chi(\wtl \D_1 ) = \chi(\wtl \D_0 ) +2 = \chi(\wtl \D ) +2 .
$$
Observe that the images of vertices $e_-', e_-''$ in $\wtl \D_1$ could be both irregular even if $e_-$ is regular in $\wtl \D$, the image of $e_+'$ in  $\wtl \D_1$ is regular and the image of $e_+''$ in  $\wtl \D_1$ is regular if and only if $e_+$ is regular in $\wtl \D$. This means that $r_0(\wtl \D_1) \le r_0(\wtl \D) +2$.
Therefore, we have
\begin{equation}\label{2e}
    \eg(\wtl \D_1) = 2- \chi(\wtl \D_1 ) +  r_0(\wtl \D_1)  \le \eg(\wtl \D) .
\end{equation}

Note that  $\wtl \D_1$ might have a \cncm\   $\wtl \D_{1,1}$  with the property that $\ph( \p G) \overset{*} = 1$ for every face $G$ in $\wtl \D_{1,1}$, i.e., $\wtl \D_1$ might lack the property (P).
Since $\wtl \D$ has property (P), it follows that there is at most one such component  $\wtl \D_{1,1}$ in $\wtl \D_1$. If $\wtl \D_{1,1}$  does exist then we take $\wtl \D_{1,1}$ out of  $\wtl \D_1$
and denote thus obtained diagram $\wtl \D_2$. If $\wtl \D_{1,1}$  does not exist then we set $\wtl \D_2 := \wtl \D_1$.
Clearly, $\wtl \D_2$ has property (P).
\smallskip

First we consider the subcase when either $\wtl \D_{1,1}$  does not exist or $\wtl \D_{1,1}$ exists and
$\chi(\wtl \D_{1,1}) \le 0$. Since $  \chi(\wtl \D_1 ) = \chi(\wtl \D_2 ) +  \chi(\wtl \D_{1,1})$, it follows from
the definitions and the inequality \eqref{2e} that $-\chi(\wtl \D_2 ) \le -\chi(\wtl \D )-2$ and
$\eg(\wtl \D_2) \le \eg(\wtl \D)$.
Hence, $\tau(\wtl \D_2 ) < \tau(\wtl \D )$ and, by the induction hypothesis, this subcase is complete.

\smallskip

Now assume that $\wtl \D_{1,1}$ exists and  $\chi(\wtl \D_{1,1}) > 0$.
Since $\wtl \D_{1,1}$ is oriented and connected, it follows that $\chi(\wtl \D_{1,1}) =2$  and $\wtl \D_{1,1}$ is a sphere.
 Since $  \chi(\wtl \D_1 ) = \chi(\wtl \D_2 ) +  \chi(\wtl \D_{1,1})$, we have $\chi(\wtl \D_2 ) = \chi(\wtl \D) $.

Let us show that $r_0(\wtl \D_2) \le r_0(\wtl \D)$. It follows from our construction that
either $r_0(\wtl \D_1) = r_0(\wtl \D)$ or  $r_0(\wtl \D_1) = r_0(\wtl \D) +2$. If
$r_0(\wtl \D_1) = r_0(\wtl \D)$ then $r_0(\wtl \D_2) \le r_0(\wtl \D)$ as desired.
Assume that  $r_0(\wtl \D_1) = r_0(\wtl \D) +2$.
Then it follows from the definitions that $r_0(\wtl \D_{1,1}) \ge 1$ because the image of the vertex
$e_-'$ in $\wtl \D_{1,1}$ is irregular. It is not difficult to show (e.g., by induction on
$(| \wtl \D_{1,1}(2) |, | \wtl \D_{1,1}(1)| )$) that the equality $r_0(\wtl \D_{1,1}) = 1$ is impossible. Therefore,
 $r_0(\wtl \D_{1,1}) \ge 2$ and we can conclude that
 $$
 r_0(\wtl \D_2) \le r_0(\wtl \D_1) -2 \le   r_0(\wtl \D) ,
 $$
 as desired.

Since $| \wtl \D_2(1) | + | \wtl \D_{1,1}(1) | = | \wtl \D(1) | $ and $| \wtl \D_{1,1}(1) | >0$, it follows from
$\chi(\wtl \D_2 ) = \chi(\wtl \D)$ and  $r_0(\wtl \D_2) \le r_0(\wtl \D)$ that $\tau(\wtl \D_2 ) < \tau(\wtl \D )$.
It is also clear that $\eg(\wtl \D_2) \le \eg(\wtl \D)$, hence, by the induction hypothesis, Case 3 is complete.

Thus in all Cases 1--3 we have been able to construct a diagram $\wtl \D'$  over $A*B$ such that $\eg(\wtl \D') \le \eg(\wtl \D)$, $\tau(\wtl \D') < \tau(\wtl \D)$ and $\wtl \D'$ contains $m$ faces $G_1, \dots, G_m$ such that
$\ph ( \p G_i) \equiv_1 u_i$, $i=1, \dots, m$.
This completes the proof of Lemma~1.
\end{proof}

\section{Car Motions}

This Section is similar to a corresponding section of \cite{FK12}
and contains necessary definitions and statements of \cite{Kl97}, \cite{Kl05} with some
simplifications.

Consider a closed map $\D$ on a closed oriented compact surface. A \emph{car} moving
around a face $F$ of $\D$  is an orientation preserving covering of the
boundary path $\p F$ of $F$ by an oriented circle $C = \mathbb R / M \mathbb Z$ called the
\emph{circle of time} and regarded as the quotient of the  real numbers $\mathbb R$ by its subgroup
 $M \mathbb Z$, where $\mathbb Z$ is the set of integers and  $M \in  \mathbb R$.

Informally, a car is a point  moving along the boundary path of a face
in counterclockwise direction (the interior of the face remains on the left)
without U-turns and stops. The motion is periodic.

The {\em degree} of a vertex $v$ of a map $\D$ is the number of oriented edges of
$\D$ whose terminal vertex is $v$. By the definition, a point in the interior of an edge of
$\D$ has {\em degree} two.

Let $v$  be a point of the
1-skeleton $\D[1]$  of  $\D$ and suppose that the number of cars being at a moment of time $t$ at the point $v$  is equal to the degree of $v$. Then $v$ is called a \emph{point of complete collision}.

A \emph{multiple car motion of period $T$} on  $\D$ is a set of  cars
$\alpha_{F,j} : C \to \p F$, defined for every face  $F$  of $\D$ and for every  $j=1,\dots,d_F$, where $d_F \ge 1$ is an integer,  such that the following hold true.

\begin{enumerate}

\item[(M1)] If $d_F >1$ then
$\alpha_{F,j}(t+T)=\alpha_{F,j+1}(t)$  for every  $t \in \mathbb R$  and
$j \in \{1,\dots,d_F\}$, here the second subscripts are modulo ${d_F}$ and addition of points
of  $C$ is defined according to $C = \mathbb R / M \mathbb Z$, where $M$ is an integer multiple of $T$.

\item[(M2)]   For every face $F$  of $\D$, there exists a partition of  $\partial F$ into $d_F$ consecutive arcs with
disjoint interiors such that, during the time interval $[0,T]$, each car
$\alpha_{F,j}$ is moving along the $j$th arc of the partition.
\end{enumerate}

\begin{lm2}[\cite{Kl97}, \cite{Kl05}]\label{lem2}
For every multiple car motion defined on a closed map $\D$ on an oriented compact surface,
the number of points of complete collision is at least
$$
\chi(\D)+\sum_{F\in \D(2)}(d_F-1),
$$
where the summation runs over all faces $F$ of $\D$.
\end{lm2}

We remark that, in articles \cite{Kl97}, \cite{Kl05},  Lemma~2 is stated and proved for
connected surfaces, but it remains valid in nonconnected case
because both parts of the inequality in Lemma~2 are additive under disjoint union.

\section{Proof of Theorem}

First we note that an arbitrary free product
$\FF = \underset {\al \in I} \ast G_\al$
of nontrivial groups $G_\al$, where $|I| >1$,
can be embedded into a free product  $A*B$ of two groups $A, B$
by means of a monomorphism $\mu : \FF \to  A*B$
in such a way that the following properties (E1)--(E2) hold true.
\begin{enumerate}
\item[(E1)]
If $w  \in \FF$ is a  reduced  word  then $\mu(w) \in  A*B$ is also reduced and
the set of finite orders of letters of  $w$
is identical to that of $\mu(w)$.

\item[(E2)]
An  element $w \in \FF$  is conjugate in $\FF$ to
an element of a free factor $G_\al$ if and only if $\mu(w)$
is conjugate in $A *B$ to an element of $A\cup B$.
\end{enumerate}

Indeed, let $A := \FF$ and let $B := F(b_\al; \al \in I)$
be a free group whose free generators are
letters $b_\al, \al \in I$.
Then the desired embedding $\mu : \FF \to  A*B$
can be defined by extending to $\FF$ the map
$\mu(g) :=  b_\al^{-1} g {b_\al}$ for every   $g \in  G_\al$.
It is easy to see that both properties (E1) and (E2) hold true.
\smallskip

Observe that if $w \in \FF$ then it follows  from property (E2) that
$$
\qp(w) \le \qp(\mu(w)) \quad \mbox{and} \quad  \mg(\mu(w)) \le \mg(w) .
$$
Hence, in view of property (E1), it suffices to prove our Theorem for
the free product $A*B$ of two factors $A, B$.

Let $w\in A*B$ be a word such that $\qp(w)$
is finite and consider a quasiperiodic  factorization for $w$ of the form
$$
w \overset { *} =  s_1 u^{n_1} s_1^{-1}   s_2 u^{n_2 }s_2^{-1} \dots  s_m u^{n_m} s_m^{-1} ,
$$
where $u$ is a cyclically reduced word, $s_j \in A*B$, $n_j >0$,
and $\qp(w) = \sum_j(n_j-1)$.
By Lemma~1, there exists a reduced
diagram $\D$  over $A*B$ such that $\D$  contains precisely
$m$ faces $F_1, \dots, F_m$ whose labels are the words
$u^{n_1}$,  $u^{n_2}, \dots, u^{n_m}$, resp.,
and
\begin{equation}\label{e2a}
    \eg(\D) \le \mg(w) .
\end{equation}

Denote
\begin{equation}\label{e3a}
 u \equiv a_1 b_1 \dots a_r b_r ,
\end{equation}
where $a_i \in A $, $b_i\in B$ and $a \ne 1$, $b \ne 1$.
\smallskip

We will now define a multiple car motion on $\D$
in the following manner. For every $j =1, \dots, m$,
there are  $n_j$ cars that move around the boundary path $\p F_j$, where $\ph(\p F_j) \equiv  u^{n_j}$,
with constant speed, one edge per  unit of time, and, at the initial moment of
time, $t=0$, the cars are located at distinct corners whose $\ph$-labels are $b_r$, here $b_r$ means the last letter of $u$, see \eqref{e3a}.
It is easy to see that this is a periodic motion with period $2r$. By Lemma~2, there are at least
$\chi(\D) + \sum_j(n_j-1)$ points of complete collision in $\D$.
\smallskip

Let us analyze where these complete collisions may occur.

First,  note that a complete collision may not occur at an interior point of an edge of $\D$.
Indeed, at every even moment of time $t=2i$, where $i \in \mathbb Z$, all cars are located at $B$-vertices,
while at every odd moment of time $t=2i+1$ all cars are located at $A$-vertices.
Therefore, during the time interval $(2i, 2i+1)$ every car is moving from a $B$-vertex to an  $A$-vertex,
while  during the time interval $(2i-1, 2i)$ every car is moving from
an $A$-vertex to a  $B$-vertex. Thus any two cars are never moving along the same edge in opposite
directions and may not collide in the interior of an edge.

Second, observe that a complete collision may not occur at a regular vertex.  To prove this claim, we note that
at every integer moment of time all
cars are located at corners with the same $\ph$-label, as denoted
in  \eqref{e3a}. More specifically,
at an even moment of time $t=2i$, where $i \in \mathbb Z$,
all cars are located at corners with $\ph$-label being
 $b_i$, as indicated in  \eqref{e3a}, here indices are modulo $r$,
and, at an odd moment of time $t=2i+1$, all
cars are located at corners with $\ph$-label being $a_i$,
as denoted in  \eqref{e3a}.
Therefore, all corners, whose vertex $v$ is
a given point of a complete collision, must have
the same $\ph$-label, as indicated in the  factorization \eqref{e3a}.
If $v$ is a regular vertex of degree $d$ then
it follows from the definition of a regular vertex that
$a_i^d \overset A = 1$ or $b_i^d \overset B = 1$ for some $i$.
Since $d$ does not exceed the number of all corners in $\D$
with $\ph$-label $a_i$ or $b_i$, it follows that
$d \le \sum_j n_j$.
However, this inequality   contradicts the  assumption that
every letter of $u$ has order greater than $\sum_j n_j$.
This contradiction proves our claim.
\smallskip

Therefore, complete collisions can only occur at irregular
vertices of $\D$. Recall that, by Lemma~2, there are at least
$\chi(\D) + \sum_j(n_j-1)$ points of complete collision in $\D$.
Hence, we conclude that the number of  irregular
vertices of $\D$ is at least $\chi(\D) + \sum_j(n_j-1)$, i.e.,
$$
r_0(\D) \ge \chi(\D) + \sum_j(n_j-1) .
$$
Therefore, it follows from the inequality \eqref{e2a}  that
$$
\qp(w) = \sum_j(n_j-1) \le -\chi(\D) + r_0(\D) = \eg(\D) -2 \le \mg(w) -2 ,
$$
as required. This completes the proof of Theorem.
\medskip

Corollaries are straightforward from the definitions and Theorem.

\medskip
{\em Acknowledgments.}   The authors thank Lvzhou Chen and Danny Calegari for helpful discussion
and remarks. The authors are grateful to the referee for useful suggestions.


\begin{thebibliography}{[10]}
\bibitem[1]{BS}
G. Baumslag and A. Steinberg,
{\em Residual nilpotence and relations in free groups},
 Bull. Amer. Math. Soc.  {\bf 70}(1964), 283--284.

\bibitem[2]{C}
L. Chen,
{\em Spectral gap of scl in free products},
Proc. Amer. Math. Soc.  {\bf 146}(2018), 3143--3151.
See also  {\tt arXiv:1611.07936  [math.GT]}.


\bibitem[3]{CG95}
A. Clifford and R. Z. Goldstein,
{\em Tesselations of $S^2$ and equations over torsion-free groups},
 Proc. Edinburgh Math. Soc. {\bf 38}(1995), 485--493.

\bibitem[4]{CG00}
A. Clifford and  R. Z. Goldstein,
{\em Equations with torsion-free coefficients},
Proc. Edinburgh Math. Soc. {\bf 43}(2000), 295--307.

\bibitem[5]{CR01}
M. M. Cohen and C. Rourke,
{\em The surjectivity problem for one-generator, one-relator extensions of
torsion-free groups},
{Geometry \& Topology} {\bf 5}(2001), 127--142.
See also arXiv:math.GR/0009101.

\bibitem[6]{CCE91}
J. A. Comerford, L. P. Comerford and C. C. Edmunds,
{\em Powers as products of commutators},
Comm. Algebra {\bf 19}(1991), 675--684.

\bibitem[7]{CER94}
L. P. Comerford,  C. C. Edmunds, and G. Rosenberger,
{\em Commutators as powers in free products of groups},
Proc. Amer. Math. Soc. {\bf 122}(1994), 47--52.

\bibitem[8]{Cull81}
M. Culler,
{\em  Using surfaces to solve equations in free groups},
Topology {\bf 20}(1981), 133--145.

\bibitem[9]{DH91}
A. J. Duncan and J. Howie,
{\em The genus problem for one-relator products of locally indicable
groups},
Math.~Z. {\bf 208}(1991), 225--237.

\bibitem[10]{FeR96}
R. Fenn and C. Rourke,
{\em Klyachko's methods and the solution of equations over torsion-free groups},
{L'Enseignment Math.} {\bf 42}(1996), 49--74.

\bibitem[11]{FeR98}
R. Fenn and C. Rourke
{\em Characterisation of a class of equations with solution over torsion-free
groups},  {``The Epstein Birthday Schrift"},
{(ed. by I. Rivin, C. Rourke and C. Series)},
{Geometry and Topology Monographs} {\bf 1}(1998), 159--166.

\bibitem[12]{FK12}
E. V. Frenkel and Ant. A. Klyachko,
{\em Commutators cannot be proper powers in metric
small-cancellation torsion-free groups}, preprint,
{\tt arXiv:1210.7908  [math.GR]}.

\bibitem[13]{FoR05}
M. Forester and C. Rourke,
Diagrams and the second homotopy group,
{\em  Comm. Anal. Geom.} {\bf 13}(2005), 801--820.
See also arXiv:math.AT/0306088.

\bibitem[14]{GK}
 R. I. Grigorchuk and P. F. Kurchanov,
{\em On the width of elements in free groups},
Ukrainian Math. J. {\bf 43}(1991), 911--918.

\bibitem[15]{How83}
J. Howie,
{\em  The solution of length three equations over groups},
{Proc. Edinburgh Math. Soc.} {\bf 26}(1983), 89--96.

\bibitem[16]{How90}
J. Howie,
{\em The quotient of a free product of groups by a single high-powered relator.
II. Fourth powers},
Proc. London Math. Soc. {\bf  61}(1990), 33--62.


\bibitem[17]{Iv}
S. V. Ivanov,
{\em The bounded and precise word problems for presentations of groups},
preprint,
{\tt arXiv:1606.08036 [math.GR]}.

\bibitem[18]{IK}
S. V. Ivanov and Ant. A. Klyachko,
{\em The asphericity and Freiheitssatz for certain LOT-presentations
of groups},
Internat. J. Algebra Comp. {\bf 11}(2001), 291--300.

\bibitem[19]{Kl93}
Ant. A. Klyachko,
A funny property of a sphere and equations over groups,
{\em  Comm. Algebra} {\bf 21}(1993), 2555--2575.

\bibitem[20]{Kl97}
Ant. A. Klyachko,
{\em Asphericity tests},
Internat. J. Algebra Comp. {\bf 7}(1997), 415--431.

\bibitem[21]{Kl05}
Ant. A. Klyachko
{\em The Kervaire--Laudenbach conjecture and presentations of simple groups},
Algebra and Logic {\bf 44}(2005), 219--242.
See also {arXiv:math.GR/0409146}.

\bibitem[22]{Kl06a}
Ant. A. Klyachko,
{\em How to generalize known results on equations over groups},
Mathematical Notes {\bf 79}(2006), 377--386.
See also
arXiv:math.GR/0406382.


\bibitem[23]{Kl06b}
Ant. A. Klyachko,
{\em The SQ-universality of one-relator relative presentations},
Sbornik Mathematics {\bf 197}(2006), 1489--1508.
See also arXiv:math.GR/0603468.

\bibitem[24]{Kl07}
Ant. A. Klyachko,
{\em Free subgroups of one-relator relative presentations},
Algebra and Logic {\bf 46}(2007), 158--162.
See also arXiv:math.GR/0510582.

\bibitem[25]{Kl09}
Ant. A. Klyachko,
{\em  The structure of one-relator relative presentations and their centres},
J. Group Theory {\bf 12}(2009), 923--947.
See also arXiv:math.GR/0701308.

\bibitem[26]{KlL12}
Ant. A. Klyachko and D. E. Lurye,
{\em  Relative hyperbolicity and similar
properties of one-generator one-relator relative presentations with
powered unimodular relator},
J. Pure Appl. Algebra {\bf 216}(2012), 524--534.
See also arXiv:1010.4220.

\bibitem[27]{Le09}
Le Thi Giang,
{\em  The relative hyperbolicity of one-relator relative presentations},
J. Group Theory  {\bf 12}(2009), 949--959.
See also arXiv:0807.2487.


\bibitem[28]{LS}
R. C. Lyndon and P. E. Schupp,
{\em Combinatorialgroup theory},
Springer-Verlag,  1977.

\bibitem[29]{Ol89}
A. Yu. Ol'shanskii,  {\em Geometry of defining relations in
groups}, Nauka, Moscow, 1989; English translation: {\em Math.
and Its Applications, Soviet series},  vol. 70, Kluwer Acad.
Publ., 1991.

\bibitem[30]{Ol}
A. Yu. Ol'shanskii, {\em On calculation of width in free groups},
London Math. Soc. Lecture Note Ser. {\bf 204}(1995), 255--258.

\bibitem[31]{Sch59}
M. P. Sch\"utzenberger, {\em  Sur  l'equation $a^{2+n}=b^{2+m}c^{2+p}$ dans
un groupe libre,}
C. R. Acad. Sci.  Paris {\bf 248}(1959), 2435--2436.
\end{thebibliography}
\end{document}